# ON THE SOLUTION OF CONSTRAINED SYLVESTER-OBSERVER EQUATION

Konstadinos H. Kiritsis

Hellenic Air Force Academy, Department of Aeronautical Sciences, Division of
Automatic Control, Dekelia Air Base, PC 13671 Acharnes, Attikis,Tatoi, Greece
E-mail: konstantinos.kyritsis@hafa.haf.gr

**Abstract**

*In this paper explicit necessary and sufficient conditions for the constrained Sylvester-observer equation are established, in order to have a solution over the field of real numbers. Furthermore, a procedure is given for the computation of the solution. Our approach is based on properties of real and polynomial matrices. Applications of the main results of this paper to linear control theory are discussed.*

## 1. INTRODUCTION

Let **R** be the field of real numbers. Furthermore, let **A** and **C** be given matrices over **R** of size *(n x n)* and *(p x n)* respectively. The problem studied in this paper can be stated as follows: Do there exist, matrices **F**, **T** and **G** over **R** of dimensions *(n-p) x (n-p)*, *(n-p) x n* and *((n-p) x p)* respectively with **F** and **T** being Hurwitz stable (i.e., all its eigenvalues have negative real parts), and of full row rank respectively such that

$$\mathbf{TA} - \mathbf{FT} = \mathbf{GC} \tag{1}$$

$$det\begin{bmatrix}\mathbf{C}\\\mathbf{T}\end{bmatrix} \neq \mathbf{0} \tag{2}$$

If so, give conditions for the existence of matrices **T**, **F** and **G** which satisfy (1) and constraint (2). The equation (1) with constraint (2) is usually called, in bibliography, "constrained Sylvester-observer equation". This term comes from the fact that equation (1) with constraint (2) arises in the design of reduced order observer for linear time-invariant systems (see Refs. [1] and [2]). Equation (1) with constraint (2) has been studied in the last years from both theoretical and computational point of view. In [3] (see also [4]) it is proven that if the pair (**A**, **C**) is observable, then there exist real matrices **T**, **F** and **G** with **T**, **F** being of full row rank and Hurwitz stable respectively which satisfy (1) and constraint (2). A criterion for choosing the matrix **F** which ensures that the solution matrix **T** of (1) is of full row rank has been presented in [5]. Furthermore an algorithm is described for the computation of solution. Elegant numerical methods have been developed for the solution of the Sylvester-observer equation in [4, 6-11]. All these methods implicitly construct a full row rank solution of (1) assuming that such a solution exists [5].

As far as we know there are no published explicit necessary and sufficient conditions for the solvability of constrained Sylvester-observer equation. This motivates the present study. In this paper, by using basic concepts and basic results from linear systems and control theory as well as of the theory of matrices, explicit necessary and sufficient conditions for the solvability of constrained Sylvester-observer equation over the field of real numbers are established. In particular, it is proved that constrained Sylvester-observer equation has a solution over the field of real numbers if and only if the pair (**A**, **C**) is detectable. Furthermore, a procedure is given for the computation of the solution.

## 2. BASIC CONCEPTS AND PRELIMINARY RESULTS

This section contains lemmas, which are needed to prove the main results of this paper and some basic concepts from linear systems and control theory as well as of the theory of matrices that are used throughout the paper. Let **R**[$s$] be the ring of polynomials with coefficients in **R**. Further let $C$ be the field of complex numbers, also let $C^+$ be the set of all complex numbers $\lambda$ with $Re(\lambda) \geq 0$. A matrix whose elements are polynomials over **R**[$s$] is termed polynomial matrix. A matrix **U**($s$) over **R**[$s$] of dimensions ($k$ x $k$) is said to be unimodular if and only if its inverse exists and is also polynomial matrix. Every polynomial matrix **W**($s$) of size ($p$ x $m$) with $rank$[**W**($s$)]=$r$, can be expressed as [12]

$$\mathbf{U}_1(s)\,\mathbf{W}(s)\,\mathbf{U}_2(s) = \mathbf{M}(s) \tag{3}$$

The polynomial matrices **U**$_1$($s$) and **U**$_2$($s$) are unimodular and the matrix **M**($s$) is given by

$$\mathbf{M}(s) = \begin{bmatrix} \mathbf{M}_r(s) & \mathbf{0} \\ \mathbf{0} & \mathbf{0} \end{bmatrix} \tag{4}$$

The non-singular polynomial matrix $\mathbf{M}_r(s)$ of size ($r$ x $r$) in (4) is given by

$$\mathbf{M}_r(s) = diag\,[a_1(s),\,a_2(s),\,....,\,a_r(s)] \tag{5}$$

The nonzero polynomials $a_i(s)$ for $i=1,2,...,\,r$ are termed invariant polynomials of $\mathbf{W}(s)$ and have the following property

$$a_i(s) \text{ divides } a_{i+1}(s),\text{ for } i=1,2,...,\,r\text{-}1 \tag{6}$$

The relationship (3) is called Smith-McMillan form of $\mathbf{W}(s)$ over $\mathbf{R}[s]$. Since the matrices $\mathbf{U}_1(s)$ and $\mathbf{U}_2(s)$ are unimodular and the polynomial matrix $\mathbf{M}_r(s)$ is non-singular, from (3) and (4) it follows that

$$rank[\mathbf{W}(s)] = rank[\mathbf{M}_r(s)] = r \tag{7}$$

Let $\mathbf{A}(s)$ and $\mathbf{B}(s)$ be matrices over $\mathbf{R}[s]$ of appropriate dimensions, if there is a matrix $\mathbf{Q}(s)$ over $\mathbf{R}[s]$ of appropriate size such that

$$\mathbf{A}(s) = \mathbf{B}(s)\mathbf{Q}(s) \tag{8}$$

Then the matrix $\mathbf{Q}(s)$ is called a right divisor of the matrix $\mathbf{A}(s)$ and the matrix $\mathbf{A}(s)$ is called a left multiple of the matrix $\mathbf{Q}(s)$ [13]. Let $\mathbf{A}(s)$ and $\mathbf{B}(s)$ be matrices over $\mathbf{R}[s]$ of appropriate dimensions, if there are matrices $\mathbf{D}(s)$, $\mathbf{A}_1(s)$ and $\mathbf{B}_1(s)$ over $\mathbf{R}[s]$ of appropriate dimensions, such that

$$\mathbf{A}(s) = \mathbf{A}_1(s)\mathbf{D}(s),\quad \mathbf{B}(s) = \mathbf{B}_1(s)\mathbf{D}(s) \tag{9}$$

Then the polynomial matrix $\mathbf{D}(s)$ is called common right divisor of polynomial matrices $\mathbf{A}(s)$ and $\mathbf{B}(s)$ [13]. A greatest common right divisor $\mathbf{D}(s)$ of two polynomial matrices $\mathbf{A}(s)$ and $\mathbf{B}(s)$ is a common right divisor which is a left multiple of every common right divisor of the matrices $\mathbf{A}(s)$ and $\mathbf{B}(s)$, [13].

Let $\mathbf{A}$ and $\mathbf{C}$ be matrices over $\mathbf{R}$ of size size *(n x n)* and *(p x n)* respectively and $\mathbf{C}$ not zero. Then there always exists a unimodular matrix $\mathbf{U}(s)$ over $\mathbf{R}[s]$ such that [13]

$$\begin{bmatrix}\mathbf{C}\\ \mathbf{I}s-\mathbf{A}\end{bmatrix} = \mathbf{U}(s)\begin{bmatrix}\mathbf{V}(s)\\ 0\end{bmatrix} \tag{10}$$

The non-singular polynomial matrix $\mathbf{V}(s)$ of size *(n x n)* is a greatest common right divisor of the polynomial matrices $\mathbf{C}$ and $[\mathbf{I}s - \mathbf{A}]$ [13]. Since the polynomial matrix $\mathbf{U}(s)$ is unimodular from (10) it follows that

$$rank\begin{bmatrix}\mathbf{C}\\ \mathbf{I}s-\mathbf{A}\end{bmatrix} = rank\begin{bmatrix}\mathbf{V}(s)\\ 0\end{bmatrix} = rank[\mathbf{V}(s)] = n \tag{11}$$

**Definition 1:** The nonzero polynomial $c(s)$ over $\mathbf{R}[s]$ is said to be strictly Hurwitz if and only if $c(s) \neq 0$, $\forall s \in C^+$.

**Definition 2:** Let $V(s)$ be a non-singular matrix over $\mathbf{R}[s]$, of size $(m \times m)$. Also let $c_i(s)$ for $i=1,2,...,m$ be the invariant polynomials of polynomial matrix $V(s)$. The polynomial matrix $V(s)$ is said to be strictly Hurwitz if and only if the polynomials $c_i(s)$ are strictly Hurwitz for every $i=1,2,...,m$, or alternatively, if and only if $det[V(s)]$ is a strictly Hurwitz polynomial.

**Definition 3:** Let $A$ and $C$ be matrices over $\mathbf{R}$ of size $(n \times n)$ and $(p \times n)$, respectively. Then the pair $(A, C)$ is said to be detectable if and only if there exists a matrix $K$ over $\mathbf{R}$ of size $(n \times p)$ such the matrix $[A+KC]$ is Hurwitz stable [14].

**Definition 4:** Let $A$ and $C$ be matrices over $\mathbf{R}$ of size $(n \times n)$ and $(p \times n)$, respectively and $C$ not zero. Then an eigenvalue $\lambda$ of the matrix $A$ is said to be observable, if and only if the following condition holds [15]:

$$rank \begin{bmatrix} C \\ I\lambda - A \end{bmatrix} = n$$

Let $A$ be a real matrix of size $(n \times n)$. The spectrum of the matrix $A$, is the set of all its eigenvalues and is denoted by $\sigma(A)$. An eigenvalue $\lambda$ of $A$ is called a stable eigenvalue if and only if $Re(\lambda) < 0$. Otherwise, the eigenvalue $\lambda$ of the matrix $A$ is said to be unstable.

**Lemma 1:** Let $A$ and $C$ be matrices over $\mathbf{R}$ of size $(n \times n)$ and $(p \times n)$, respectively and $C$ not zero. Further, let $\sigma(A)$ be the spectrum of the matrix $A$. The pair $(A, C)$ is detectable, if and only if one of the following equivalent conditions holds [16]:

(a) $rank \begin{bmatrix} C \\ Is - A \end{bmatrix} = n$, $\forall s \in C^+$

(b) $rank \begin{bmatrix} C \\ I\lambda - A \end{bmatrix} = n$, $\forall \lambda \in \sigma(A)$ with $Re(\lambda) \geq 0$.

From condition (b) of Lemma 1 it follows that the pair $(A, C)$ is detectable if and only if all unstable eigenvalues of the matrix $A$ are observable, [16]. The following lemma is taken from [17].

**Lemma 2:** Let $V(s)$ be a non-singular polynomial matrix over $\mathbf{R}[s]$, of size $(m \times m)$. Also let $c_i(s)$ for $i=1,2,...,m$ be the invariant polynomials of the polynomial matrix $V(s)$. The polynomial matrix $V(s)$ is strictly Hurwitz if and only if the following condition holds

(a) $rank[V(s)] = m$, $\forall s \in C^+$

***Proof:*** Let $\mathbf{V}(s)$ be a strictly Hurwitz polynomial matrix of size $(m \times m)$ with invariant polynomials $c_i(s)$ for $i=1,2,...,m$. From Definition 2 it follows that the polynomials $c_i(s)$ are strictly Hurwitz for every $i=1,2,...,m$ and therefore from Definition 1 it follows that

$$c_i(s) \neq 0, \forall s \in C^+, \forall i=1,2,...,m \tag{12}$$

We define the polynomial matrix

$$\mathbf{V}_m(s) = \text{diag }[c_1(s), c_2(s), ....c_m(s)] \tag{13}$$

From (12) and (13) it follows that

$$rank[\mathbf{V}_m(s)] = rank\{diag\,[c_1(s),c_2(s),....,c_m(s)]\}=m, \forall s \in C^+ \tag{14}$$

The Smith-McMillan form of polynomial matrix $\mathbf{V}(s)$ over $\mathbf{R}[s]$ is given by

$$\mathbf{K}(s)\,\mathbf{V}(s)\,\mathbf{L}(s) = \mathbf{V}_m(s) \tag{15}$$

where $\mathbf{K}(s)$ and $\mathbf{L}(s)$ are unimodular matrices. Since the matrices $\mathbf{K}(s)$ and $\mathbf{L}(s)$ are unimodular, from (7), (14) and (15) we have that

$$rank[\mathbf{V}(s)]=m, \forall s \in C^+ \tag{16}$$

This is condition (a) of the Lemma. To prove sufficiency, we assume that condition (a) holds. Using (7) from (13) and (15) we have that

$$rank[\mathbf{V}(s)] = rank[\mathbf{V}_m(s)] = rank\{diag[c_1(s),c_2(s),....,c_m(s)]\} = m \tag{17}$$

Since by assumption condition (a) holds we have that

$$rank[\mathbf{V}(s)]=m, \forall s \in C^+ \tag{18}$$

Relationships (17) and (18) imply

$$rank[\mathbf{V}_m(s)] = rank\{diag\,[c_1(s),c_2(s),....,c_m(s)]\}=m, \forall s \in C^+ \tag{19}$$

From (19) it follows that

$$c_i(s) \neq 0, \forall s \in C^+, \forall i=1,2,...,m \tag{20}$$

Relationship (20) and Definition 1 imply that polynomials $c_i(s)$ are strictly Hurwitz for every $i=1,2,...,m$, and therefore according to Definition 2 the non-singular polynomial matrix $\mathbf{V}(s)$ over $\mathbf{R}[s]$, is strictly Hurwitz. This completes the proof.

The following lemma is taken from [17].

**Lemma 3:** Let $\mathbf{A}$ and $\mathbf{C}$ be matrices over $\mathbf{R}$ matrices of size $(n \times n)$, $(p \times n)$, respectively and $\mathbf{C}$ not zero. Further let $\mathbf{V}(s)$ be a greatest common right divisor of polynomial matrices $[\mathbf{I}s - \mathbf{A}]$ and $\mathbf{C}$ of size $(n \times n)$. The pair $(\mathbf{A}, \mathbf{C})$ is detectable if and only if the following condition holds:

(a) The polynomial matrix $\mathbf{V}(s)$ is strictly Hurwitz.

*Proof*: Let the pair (**A**, **C**) is detectable. Then from condition (a) of Lemma 1 it follows that

$$rank\begin{bmatrix} \mathbf{C} \\ \mathbf{I}s - \mathbf{A} \end{bmatrix} = n, \quad \forall s \in C^+ \tag{21}$$

Since the polynomial matrix **V**(s) is the greatest common right divisor of polynomial matrices [**I**s − **A**] and **C**, from (10) it follows that there exists a unimodular matrix **U**(s) such that

$$\begin{bmatrix} \mathbf{C} \\ \mathbf{I}s - \mathbf{A} \end{bmatrix} = \mathbf{U}(s)\begin{bmatrix} \mathbf{V}(s) \\ \mathbf{0} \end{bmatrix} \tag{22}$$

Since the polynomial matrix **U**(s) is unimodular from (11) and (22) it follows that

$$rank\begin{bmatrix} \mathbf{C} \\ \mathbf{I}s - \mathbf{A} \end{bmatrix} = rank\begin{bmatrix} \mathbf{V}(s) \\ \mathbf{0} \end{bmatrix} = rank[\mathbf{V}(s)] \tag{23}$$

From relationships (21) and (23) we have that

$$rank[\mathbf{V}(s)] = n, \quad \forall s \in C^+ \tag{24}$$

Relationship (24) and condition (a) of Lemma 2 imply that the polynomial matrix **V**(s) is strictly Hurwitz. This is condition (a) of the Lemma. To prove sufficiency, we assume that the polynomial matrix **V**(s) is strictly Hurwitz. Then from Lemma 2 it follows that

$$rank[\mathbf{V}(s)] = n, \quad \forall s \in C^+ \tag{25}$$

Since the polynomial matrix **V**(s) is the greatest common right divisor of polynomial matrices [**I**s − **A**] and **C**, from (23) and (25) it follows that

$$rank\begin{bmatrix} \mathbf{C} \\ \mathbf{I}s - \mathbf{A} \end{bmatrix} = n, \quad \forall s \in C^+ \tag{26}$$

Condition (a) of Lemma 1 and (26) imply that the pair (**A**, **C**) is detectable. This completes the proof.

Let **A** and **C** be matrices over **R** matrices of size *(n x n)*, *(p x n)*, respectively with *rank*[**C**]=*p*. Then there exists a non-singular matrix **L** over **R** of size *(n x n)* such that

$$\mathbf{CL} = \mathbf{C}_1 \tag{27}$$

$$\mathbf{L}^{-1}\mathbf{AL} = \mathbf{A}_1 \tag{28}$$

The matrices $\mathbf{C}_1$ and $\mathbf{A}_1$ are given by

$$\mathbf{C}_1 = [\mathbf{I}_p, \mathbf{0}] \tag{29}$$

$$\mathbf{A}_1 = \begin{bmatrix} \mathbf{A}_{11} & \mathbf{A}_{12} \\ \mathbf{A}_{21} & \mathbf{A}_{22} \end{bmatrix} \tag{30}$$

where $I_p$ is the identity matrix of size *(p x p)*, and $A_{11}$, $A_{12}$, $A_{21}$ and $A_{22}$ are matrices over **R** of dimensions *(p x p)*, *(p x (n −p))*, *((n-p) x p)* and *((n-p) x (n-p))* respectively.

**Lemma 4:** Let **A** and **C** be matrices over **R** matrices of size *(n x n)*, *(p x n)* respectively with *rank*[**C**]=*p*. The pair (**A**, **C**) is detectable if and only if the following condition holds:

(a) The pair $(A_{22}, A_{12})$ is detectable.

*Proof*: Let the pair (**A**, **C**) is detectable. Then from condition (a) of Lemma 1 it follows that

$$rank \begin{bmatrix} C \\ Is - A \end{bmatrix} = n, \ \forall s \in C^+ \tag{31}$$

From (27), (28), (29) and (30) we have that

$$C = [I_p, 0] L^{-1} = C_1 L^{-1}, \ A = L A_1 L^{-1} = L \begin{bmatrix} A_{11} & A_{12} \\ A_{21} & A_{22} \end{bmatrix} L^{-1} \tag{32}$$

Using (32) and after simple algebraic manipulations the (31) can be expressed as follows

$$rank \begin{bmatrix} C \\ Is - A \end{bmatrix} = rank[diag[I_p, L] \begin{bmatrix} C_1 \\ Is - A_1 \end{bmatrix} L^{-1}] =$$

$$= rank[diag[I_p, L] \begin{bmatrix} I_p & 0 \\ I_p s - A_{11} & -A_{12} \\ -A_{21} & I_{n-p} s - A_{22} \end{bmatrix} L^{-1}] = n, \ \forall s \in C^+ \tag{33}$$

Since the matrices, $diag[I_p, L]$ and $L^{-1}$ are non-singular, from (33) it follows that

$$rank \begin{bmatrix} I_p & 0 \\ I_p s - A_{11} & -A_{12} \\ -A_{21} & I_{n-p} s - A_{22} \end{bmatrix} = n, \ \forall s \in C^+ \tag{34}$$

Since the *n* columns of the matrix in the left side of (34) are linearly independent over C, $\forall s \in C^+$, a subset of these columns consisting of the last *(n-p)* columns must be also linearly independent over C, $\forall s \in C^+$; therefore

$$rank \begin{bmatrix} 0 \\ -A_{12} \\ I_{n-p} s - A_{22} \end{bmatrix} = rank \begin{bmatrix} -A_{12} \\ I_{n-p} s - A_{22} \end{bmatrix} = rank \begin{bmatrix} A_{12} \\ I_{n-p} s - A_{22} \end{bmatrix} = (n - p), \forall s \in C^+ \tag{35}$$

Relationship (35) and condition (a) of Lemma 1 imply that the pair $(A_{22}, A_{12})$ is detectable. This is condition (a). To prove sufficiency, we assume that the pair $(A_{22}, A_{12})$ is detectable. Let $\Gamma(s)$ be a greatest common right divisor of polynomial matrices

$[I_{n-p}s - A_{22}]$ and $A_{12}$. Then from (10) it follows that there exists a unimodular matrix $U(s)$ of size $(n \times n)$ such that

$$\begin{bmatrix} -A_{12} \\ I_{n-p}s - A_{22} \end{bmatrix} = U(s) \begin{bmatrix} \Gamma(s) \\ 0 \end{bmatrix} \tag{36}$$

We define the following matrices

$$\begin{bmatrix} I_p s - A_{11} \\ -A_{21} \end{bmatrix} = U(s) \begin{bmatrix} E(s) \\ Z(s) \end{bmatrix} \tag{37}$$

$$\Delta(s) = \begin{bmatrix} I_p & 0 \\ E(s) & I_{n-p} \end{bmatrix} \tag{38}$$

$$H(s) = [Z(s), 0] \tag{39}$$

$$M(s) = diag[I_p, \Gamma(s)] \tag{40}$$

where $\Delta(s)$ is a unimodular matrix of size $(n \times n)$ and $E(s)$, $Z(s)$, $H(s)$ and $M(s)$ are polynomial matrices of size $((n-p) \times p)$, $(p \times p)$, $(p \times n)$ and $(n \times n)$ respectively. Using (29), (30), (37), (38), (39) and (40) we have that

$$\begin{bmatrix} C_1 \\ Is - A_1 \end{bmatrix} = \begin{bmatrix} I_p & 0 \\ I_p s - A_{11} & -A_{12} \\ -A_{21} & I_{n-p}s - A_{22} \end{bmatrix} = diag[I_p, U(s)] \begin{bmatrix} I_p & 0 \\ E(s) & \Gamma(s) \\ Z(s) & 0 \end{bmatrix} =$$

$$= diag[I_p, U(s)] \, diag[\Delta(s), I_p] \begin{bmatrix} I_p & 0 \\ 0 & I_{n-p} \\ Z(s) & 0 \end{bmatrix} diag[I_p, \Gamma(s)] =$$

$$= diag[I_p, U(s)] \, diag[\Delta(s), I_p] \begin{bmatrix} I_n \\ H(s) \end{bmatrix} diag[I_p, \Gamma(s)] =$$

$$= diag[I_p, U(s)] \, diag[\Delta(s), I_p] \begin{bmatrix} I_n & 0 \\ H(s) & I_p \end{bmatrix} \begin{bmatrix} I_n \\ 0 \end{bmatrix} M(s) \tag{41}$$

Since the matrices $diag[I_p, U(s)]$, $diag[\Delta(s), I_p]$ and $\begin{bmatrix} I_n & 0 \\ H(s) & I_p \end{bmatrix}$ are all unimodular, their product

$$diag[I_p, U(s)] \, diag[\Delta(s), I_p] \begin{bmatrix} I_n & 0 \\ H(s) & I_p \end{bmatrix}$$

is also a unimodular matrix and therefore from (10) it follows that the matrix $M(s)$ in (41) is a greatest common right divisor of the polynomial matrices $C_1$, $[Is - A_1]$. Detectability of the pair $(A_{22}, A_{12})$ and Lemma 3 imply that the matrix $\Gamma(s)$ in (36) is strictly Hurwitz. From (40) we have that

$$det[M(s)] = det[\Gamma(s)] \tag{42}$$

Since the matrix $\boldsymbol{\Gamma}(s)$ is strictly Hurwitz, from Definition 2 it follows that the $det[\boldsymbol{\Gamma}(s)]$ is strictly Hurwitz polynomial and therefore from (42) and Definition 2 it follows that the matrix $\mathbf{M}(s)$ is strictly Hurwitz. Since $\mathbf{M}(s)$ is strictly Hurwitz, from Lemma 3 it follows that the pair $(\mathbf{A}_1, \mathbf{C}_1)$ is detectable. Since the pair $(\mathbf{A}_1, \mathbf{C}_1)$ is detectable and the detectability is invariant under similarity transformation [16], from (32) it follows that the pair (**A**, **C**) is also detectable. This completes the proof.

**Lemma 5.** Let **A** and **C** be matrices over **R** of size *(n x n)* and *(p x n)*, respectively. Then the pair (**A**, **C**) is observable if and only if for every monic polynomial $c(s)$ over **R**[*s*] of degree *n* there exists a matrix **K** over **R** of size *( n x p)*, such that the matrix [**A**+**KC**] has characteristic polynomial $c(s)$ [12].

The standard decomposition of unobservable systems given in the following Lemma was first published by Kalman in [18] and can also be found in any standard book of linear systems theory.

**Lemma 6:** Let **A** and **C** be matrices over **R** of size *(n x n)* and *(p x n)*, respectively. Further, let the pair (**A**, **C**) is unobservable and **C** not zero. Then there exists a non-singular matrix **T** of size *(n x n)* such that

$$\mathbf{T}^{-1}\mathbf{A}\mathbf{T} = \begin{bmatrix} \mathbf{A}_{11} & \mathbf{0} \\ \mathbf{A}_{21} & \mathbf{A}_{22} \end{bmatrix},$$

$$\mathbf{CT} = [\,\mathbf{C}_1,\,\mathbf{0}\,]$$

The pair $(\mathbf{A}_{11}, \mathbf{C}_1)$ is observable and the eigenvalues of the matrix $\mathbf{A}_{22}$ are the unobservable eigenvalues of the pair (**A**, **C**).

**Lemma 7:** Let **A** and **C** be matrices over **R** of size *(n x n)*, *(p x n)*, respectively and **C** not zero. Further let

$$\mathbf{A} = \mathbf{T}\begin{bmatrix} \mathbf{A}_{11} & \mathbf{0} \\ \mathbf{A}_{21} & \mathbf{A}_{22} \end{bmatrix}\mathbf{T}^{-1}, \quad \mathbf{C} = [\,\mathbf{C}_1,\,\mathbf{0}\,]\,\mathbf{T}^{-1}$$

with $(\mathbf{A}_{11}, \mathbf{C}_1)$ observable. If the pair (**A**, **C**) is detectable then the matrix $\mathbf{A}_{22}$ is Hurwitz stable [16].

The following lemma and its proof are based on the results of [12].

**Lemma 8:** Let **A** and **C** be matrices over **R** of size *(n x n)*, *(p x n)*, respectively and **C** not zero. Further, let the pair (**A**, **C**) be detectable. Then there exists a matrix **K** over **R** of size *(n x p)* such that the matrix [**A**+**KC**] is Hurwitz stable.

***Proof***: Let the pair (**A**, **C**) is detectable. Detectability of the pair (**A**, **C**) implies that the pair (**A**, **C**) is either observable or unobservable with stable unobservable eigenvalues. If the pair (**A**, **C**) is observable, then from Lemma 5 it follows that there exists a matrix **K** of appropriate size over **R** such that

$$det[\mathbf{I}s - \mathbf{A} - \mathbf{KC}] = c(s)$$

where $c(s)$ be an arbitrary monic, strictly Hurwitz polynomial over **R**[$s$] of degree $n$. Since $c(s)$ is a strictly Hurwitz polynomial over **R**[$s$] from Definition 3 it follows that matrix [**A**+**KC**] is Hurwitz stable. Since the notion of observability is a dual of controllability (i.e., observability of the pair (**A**, **C**) implies controllability of the pair(**A**$^T$, **C**$^T$)) [12], the matrix **K** can be calculated using known methods for the solution of pole assignment problem by state feedback [12].

If the pair (**A**, **C**) is unobservable with stable unobservable eigenvalues, then from Lemma 6 and Lemma 7 it follows that there exists a matrix **T** such that

$$\mathbf{A} = \mathbf{T}\begin{bmatrix} \mathbf{A}_{11} & \mathbf{0} \\ \mathbf{A}_{21} & \mathbf{A}_{22} \end{bmatrix}\mathbf{T}^{-1}, \quad \mathbf{C} = [\mathbf{C}_1, \mathbf{0}]\,\mathbf{T}^{-1}$$

The pair (**A**$_{11}$, **C**$_1$) is observable and the matrix **A**$_{22}$ is Hurwitz stable. Hurwitz stability of the matrix **A**$_{22}$ and Definition 3 imply that the polynomial $\chi(s)$ given by

$$det[\mathbf{I}s - \mathbf{A}_{22}] = \chi(s)$$

is a strictly Hurwitz polynomial. Observability of the pair (**A**$_{11}$, **C**$_1$) and Lemma 5 imply the existence of a matrix **K**$_1$ over **R** of appropriate dimensions such that

$$det[\mathbf{I}s - \mathbf{A}_{11} - \mathbf{K}_1\mathbf{C}_1] = \varphi(s)$$

where $\varphi(s)$ is an arbitrary monic, strictly Hurwitz polynomial over **R**[$s$] of appropriate degree. Since the notion of observability is a dual of controllability (i.e., observability of the pair (**A**$_{11}$, **C**$_1$) implies controllability of the pair(**A**$_{11}^T$, **C**$_1^T$) ), the matrix **K**$_1$ can be calculated using known methods for the solution of pole assignment problem by state feedback [12]. Putting

$$\mathbf{K} = \mathbf{T}\begin{bmatrix} \mathbf{K}_1 \\ \mathbf{0} \end{bmatrix}$$

we obtain

$$\mathbf{A} + \mathbf{KC} = \mathbf{T}\begin{bmatrix} \mathbf{A}_{11} & \mathbf{0} \\ \mathbf{A}_{21} & \mathbf{A}_{22} \end{bmatrix}\mathbf{T}^{-1} + \mathbf{T}\begin{bmatrix} \mathbf{K}_1 \\ \mathbf{0} \end{bmatrix}[\mathbf{C}_1, \mathbf{0}]\,\mathbf{T}^{-1} =$$

$$= \mathbf{T}\begin{bmatrix} \mathbf{A}_{11} + \mathbf{K}_1\mathbf{C}_1 & \mathbf{0} \\ \mathbf{A}_{21} & \mathbf{A}_{22} \end{bmatrix}\mathbf{T}^{-1}$$

From the last relationships it follows that

$$det[\mathbf{I}s-\mathbf{A}-\mathbf{KC}] = \varphi(s)\chi(s)$$

Since the polynomials $\chi(s)$ and $\varphi(s)$ are strictly Hurwitz, the polynomial $\varphi(s)\chi(s)$ is also a strictly Hurwitz polynomial and therefore, from Definition 3, it follows that matrix $[\mathbf{A}+\mathbf{KC}]$ is Hurwitz stable. This completes the proof.

## 3. MAIN RESULTS

The theorem that follows is the main result of this paper and gives the necessary and sufficient conditions for the constrained Sylvester-observer equation in order, to have a solution over the field of real numbers. Without any loss of generality in what follows we assume that $rank[\mathbf{C}] = p$.

**Theorem 1.** There exists matrices **F**, **T** and **G** with **F**, **T** being Hurwitz stable and of full row rank respectively which satisfy equation (1) and constraint (2), if and only if the following condition holds:

(a) The pair (**A**, **C**) is detectable.

*Proof*: Let there exists matrices **F**, **T** and **G** with **F**, **T** being Hurwitz stable and of full row rank respectively which satisfy equation (1) and constraint (2). Also let **L** a non-singular matrix over **R** of size *(n x n)* as in (27), (28), (29) and (30). From (27), (28), (29) and (30) we have that

$$\mathbf{C} = [\mathbf{I}_p , \mathbf{0}] \mathbf{L}^{-1} \quad (43)$$

$$\mathbf{A} = \mathbf{L} \begin{bmatrix} \mathbf{A}_{11} & \mathbf{A}_{12} \\ \mathbf{A}_{21} & \mathbf{A}_{22} \end{bmatrix} \mathbf{L}^{-1} \quad (44)$$

Let

$$\mathbf{T} = [\mathbf{T}_1, \mathbf{T}_2] \mathbf{L}^{-1} \quad (45)$$

where $\mathbf{T}_1$ and $\mathbf{T}_2$ are real matrices of size $((n\text{-}p) \times p)$ and $((n\text{-}p) \times (n\text{-}p))$ respectively. Using (43) and (45) the matrix $\begin{bmatrix} \mathbf{C} \\ \mathbf{T} \end{bmatrix}$ can be expressed as follows

$$\begin{bmatrix} \mathbf{C} \\ \mathbf{T} \end{bmatrix} = \begin{bmatrix} \mathbf{I}_p & \mathbf{0} \\ \mathbf{T}_1 & \mathbf{T}_2 \end{bmatrix} \mathbf{L}^{-1} \quad (46)$$

Since according to (2) the matrix $\begin{bmatrix} \mathbf{C} \\ \mathbf{T} \end{bmatrix}$ is non-singular, from (46) it follows that the matrix $\mathbf{T}_2$ is non-singular; therefore

$$det[\mathbf{T}_2] \neq 0 \quad (47)$$

Substituting (43), (44) and (45) to (1) we obtain:

$$[\mathbf{T}_1, \mathbf{T}_2] \mathbf{L}^{-1} \mathbf{L} \begin{bmatrix} \mathbf{A}_{11} & \mathbf{A}_{12} \\ \mathbf{A}_{21} & \mathbf{A}_{22} \end{bmatrix} \mathbf{L}^{-1} - \mathbf{F}[\mathbf{T}_1, \mathbf{T}_2] \mathbf{L}^{-1} = \mathbf{G}[\mathbf{I}_p, \mathbf{0}] \mathbf{L}^{-1} \quad (48)$$

Relationship (48) can be rewritten as follows

$$[\mathbf{T}_1, \mathbf{T}_2] \begin{bmatrix} \mathbf{A}_{11} & \mathbf{A}_{12} \\ \mathbf{A}_{21} & \mathbf{A}_{22} \end{bmatrix} - \mathbf{F}[\mathbf{T}_1, \mathbf{T}_2] = \mathbf{G}[\mathbf{I}_p, \mathbf{0}] \quad (49)$$

From (49) it follows that

$$\mathbf{G} = \mathbf{T}_1 \mathbf{A}_{11} + \mathbf{T}_2 \mathbf{A}_{21} - \mathbf{F}\mathbf{T}_1 \quad (50)$$

$$\mathbf{F}\mathbf{T}_2 = \mathbf{T}_2 \mathbf{A}_{22} + \mathbf{T}_1 \mathbf{A}_{12} \quad (51)$$

Condition (47) implies that the matrix $\mathbf{T}_2$ is invertible. Using this fact (51) can be rewritten as follows

$$\mathbf{T}_2^{-1}\mathbf{F}\mathbf{T}_2 = \mathbf{A}_{22} + \mathbf{T}_2^{-1}\mathbf{T}_1 \mathbf{A}_{12} \quad (52)$$

Since by assumption the matrix $\mathbf{F}$ is Hurwitz stable, the matrix $\mathbf{T}_2^{-1}\mathbf{F}\mathbf{T}_2$ must be also Hurwitz stable. Hurwitz stability of $\mathbf{T}_2^{-1}\mathbf{F}\mathbf{T}_2$ implies Hurwitz stability of $[\mathbf{A}_{22} + \mathbf{T}_2^{-1}\mathbf{T}_1\mathbf{A}_{12}]$. Definition 3 and Hurwitz stability of the matrix $[\mathbf{A}_{22} + \mathbf{T}_2^{-1}\mathbf{T}_1\mathbf{A}_{12}]$ imply detectability of the pair $(\mathbf{A}_{22}, \mathbf{A}_{12})$. Detectability of the pair $(\mathbf{A}_{22}, \mathbf{A}_{12})$ and Lemma 4 imply detectability of the pair $(\mathbf{A}, \mathbf{C})$. This is condition (a) of the Theorem. To prove sufficiency, we assume that condition (a) holds. Detectability of the pair $(\mathbf{A}, \mathbf{C})$ and Lemma 4 imply detectability of the pair $(\mathbf{A}_{22}, \mathbf{A}_{12})$. Since the the pair $(\mathbf{A}_{22}, \mathbf{A}_{12})$ is detectable, from Lemma 8 it follows that there exists a real matrix $\mathbf{K}$ of appropriate dimensions such that the matrix $[\mathbf{A}_{22} + \mathbf{K}\mathbf{A}_{12}]$ is Hurwitz stable, that is

$$det[\mathbf{I}_{n-p}s - \mathbf{A}_{22} - \mathbf{K}\mathbf{A}_{12}] = c(s) \quad (53)$$

where $c(s)$ is a strictly Hurwitz polynomial over $\mathbf{R}[s]$ of degree $(n-p)$. The matrix $\mathbf{K}$ can be calculated as in the proof of Lemma 8. Let

$$\mathbf{T} = [\mathbf{T}_1, \mathbf{T}_2]\mathbf{L}^{-1} = [\mathbf{K}, \mathbf{I}_{n-p}] \mathbf{L}^{-1} \quad (54)$$

$$\mathbf{F} = [\mathbf{A}_{22} + \mathbf{K}\mathbf{A}_{12}] \quad (55)$$

$$\mathbf{G} = \mathbf{K}\mathbf{A}_{11} + \mathbf{A}_{21} - [\mathbf{A}_{22} + \mathbf{K}\mathbf{A}_{12}]\mathbf{K} \quad (56)$$

Using (54), (55), (56), (43) and (44) we have that

$$\mathbf{TA} - \mathbf{FT} = [\mathbf{K}, \mathbf{I}_{n-p}]\mathbf{L}^{-1}\mathbf{A} - [\mathbf{A}_{22} + \mathbf{K}\mathbf{A}_{12}][\mathbf{K}, \mathbf{I}_{n-p}]\mathbf{L}^{-1} =$$

$$= \{[\mathbf{K}, \mathbf{I}_{n-p}]\mathbf{L}^{-1}\mathbf{L}\begin{bmatrix} \mathbf{A}_{11} & \mathbf{A}_{12} \\ \mathbf{A}_{21} & \mathbf{A}_{22} \end{bmatrix}\mathbf{L}^{-1} - [\mathbf{A}_{22} + \mathbf{K}\mathbf{A}_{12}][\mathbf{K}, \mathbf{I}_{n-p}]\}\mathbf{L}^{-1} =$$

$$= \{[\mathbf{K}, \mathbf{I}_{n-p}] \begin{bmatrix} \mathbf{A}_{11} & \mathbf{A}_{12} \\ \mathbf{A}_{21} & \mathbf{A}_{22} \end{bmatrix} - [\mathbf{A}_{22} + \mathbf{K}\mathbf{A}_{12}][\mathbf{K}, \mathbf{I}_{n-p}] \} \mathbf{L}^{-1} =$$

$$= \{[(\mathbf{K}\mathbf{A}_{11} + \mathbf{A}_{21}), (\mathbf{A}_{22} + \mathbf{K}\mathbf{A}_{12})] - [(\mathbf{A}_{22} + \mathbf{K}\mathbf{A}_{12})\mathbf{K}, (\mathbf{A}_{22} + \mathbf{K}\mathbf{A}_{12})]\} \mathbf{L}^{-1} =$$

$$= \{[(\mathbf{K}\mathbf{A}_{11} + \mathbf{A}_{21}) - (\mathbf{A}_{22} + \mathbf{K}\mathbf{A}_{12})\mathbf{K}, \mathbf{0}]\}\mathbf{L}^{-1} =$$

$$= \{[(\mathbf{K}\mathbf{A}_{11} + \mathbf{A}_{21}) - (\mathbf{A}_{22} + \mathbf{K}\mathbf{A}_{12})\mathbf{K}][\mathbf{I}_p, \mathbf{0}]\} \mathbf{L}^{-1} = \mathbf{GC} \tag{57}$$

Using (43) and (54) we have that

$$\begin{bmatrix} \mathbf{C} \\ \mathbf{T} \end{bmatrix} = \begin{bmatrix} \mathbf{I}_p & \mathbf{0} \\ \mathbf{K} & \mathbf{I}_{n-p} \end{bmatrix} \mathbf{L}^{-1} \tag{58}$$

From (54) it follows that the matrix **T** is of full row rank and from (53) and (55) we have that the matrix **F** is Hurwitz stable. Relationships (57) and (58) imply that the matrices **T**, **F** and **G** given by (54), (55) and (56) satisfy equation (1) and constraint (2). This completes the proof.

The sufficiency part of the proof of Theorem 1 provides a construction of the matrices **T**, **F** and **G** which satisfy equation (1) and constrained (2). The major steps of this construction are given below.

**CONSTRUCTION**

*Given*: **A** and **C** with *rank*[**C**] = *p*.

*Find*: **T**, **F** and **G**

*Step 1:* Check condition (a) of Theorem 1. If this condition is satisfied go to *Step 2*. If condition (a) is not satisfied, the solution of equation (1) with constrained (2) is impossible.

*Step 2:* Find a non-singular matrix **L** over **R** of dimensions (*n x n*) such that

$$\mathbf{CL} = [\mathbf{I}_p, \mathbf{0}]$$

*Step 3:* Calculate the following real matrix

$$\mathbf{L}^{-1}\mathbf{AL} = \begin{bmatrix} \mathbf{A}_{11} & \mathbf{A}_{12} \\ \mathbf{A}_{21} & \mathbf{A}_{22} \end{bmatrix},$$

*Step 4:* Detectability of the pair (**A**, **C**) and Lemma 4 imply detectability of the pair ($\mathbf{A}_{22}$, $\mathbf{A}_{12}$). Since the the pair ($\mathbf{A}_{22}$, $\mathbf{A}_{12}$) is detectable from Lemma 8 it follows that there exists a real matrix **K** of appropriate size such that the matrix

$$[\mathbf{A}_{22} + \mathbf{K}\mathbf{A}_{12}]$$

is Hurwitz stable. The matrix **K** can be calculated as in the proof of Lemma 8.

*Step 5:* Put

$$T = [K, I_{n-p}] \, L^{-1}$$

$$F = [A_{22} + KA_{12}]$$

$$G = KA_{11} + A_{21} - [A_{22} + KA_{12}]K$$

## 4. APPLICATION TO CONTROL THEORY

Consider a linear time-invariant system described by the following state-space equations

$$\dot{x}(t) = Ax(t) + Bu(t) \tag{59a}$$

$$y(t) = Cx(t) \tag{59b}$$

where **A**, **B** and **C** are real matrices of size *(n x n), (n x m)* and *(p x n)* respectively, **x**(*t*) is the state vector of dimensions *(n x 1)*, **u**(*t*) is the vector of inputs of dimensions *(m x1)* and **y**(*t*) the vector of outputs of dimensions *(p x1)*. Consider also a linear time-invariant system described by the following state-space equations [1]

$$\dot{z}(t) = Fz(t) + Gy(t) + Pu(t) \tag{60}$$

where **F**, **G** and **P** are real matrices of size *((n-p) x (n-p))*, *((n-p) x p)* and *((n-p) x m)* respectively and **z**(*t*) is the state vector of dimensions *((n-p) x 1)*.

**Theorem 2.** The system (60) is an observer of order *(n-p)* of the system described by equations (59a) and (59b) if the following condition holds:

(a) The pair (**A**, **C**) is detectable.

*Proof*: Detectability of the pair (**A**, **C**) and Theorem 1 imply the existence of a Hurwitz stable matrix **F**, of full row rank matrix **T** and the matrix **G** of dimensions *((n-p) x (n-p))*, *((n-p) x n)* and *((n-p) x p)* respectively such that

$$TA - FT = GC \tag{61}$$

$$det\begin{bmatrix}C\\T\end{bmatrix} \neq 0 \tag{62}$$

Let the state estimation error be defined as follows [1]

$$e(t) = z(t) - Tx(t) \tag{63}$$

Then by taking the derivative of (63) and using (59a) and (60) we have that

$$\dot{e}(t) = \dot{z}(t) - T\dot{x}(t) = Fz(t) + Gy(t) + Pu(t) - T[Ax(t) + Bu(t)] \tag{64}$$

Substituting (59b) while adding and subtracting **FT**x(t) in (64) [1], we have that

$$\dot{e}(t) = Fe(t) - [TA - FT - GC]x(t) + [P - TB]u(t) \tag{65}$$

Let
$$\mathbf{P} = \mathbf{TB} \quad (66)$$

Using (66) and (61), from (65) we have that

$$\dot{\mathbf{e}}(t) = \mathbf{Fe}(t) \quad (67)$$

Since the matrix **F** is Hurwitz stable, from (67) we have that

$$lim_{t \to +\infty} \mathbf{e}(t) = \mathbf{0} \quad (68)$$

for any **x**(0), **z**(0), **e**(0) and **u**(t). Hence, from (63) and (68) it follows that **z**(t) is an estimate of **Tx**(t) and therefore the system (60) is an observer of order (*n-p*) of the system described by equations (59a) and (59b) [1]. The observer matrices **F**, **G** and **P** can be calculated by solving equation (61) and using relationship (66). Taking into account eqns. (62), (63), (68) and (59b), one can always find an estimate $\hat{\mathbf{x}}(t)$ of **x**(t) of the system described by equations (59a) and (59b) as follows [1]

$$\hat{\mathbf{x}}(t) = \begin{bmatrix} \mathbf{C} \\ \mathbf{T} \end{bmatrix}^{-1} \begin{bmatrix} \mathbf{y}(t) \\ \mathbf{z}(t) \end{bmatrix} \quad (69)$$

This completes the proof.

## 5. CONCLUSIONS

In this paper explicit necessary and sufficient conditions for the constrained Sylvester-observer equation in order, to have a real solution, are established. The proof of the main results of this paper is constructive and furnishes a procedure for the computation of solution. The results of this paper are useful in studying the problems of linear control theory that can be converted to the solution of constrained Sylvester-observer equation. A typical example is the problem of the existence and construction of reduced order observer for the state estimation of continuous-time linear time invariant systems.